\expandafter\ifx\csname amssym.def\endcsname\relax \else\endinput\fi
%
%  Store the catcode of the @ in the csname so that it can be restored later.
\expandafter\edef\csname amssym.def\endcsname{%
       \catcode`\noexpand\@=\the\catcode`\@\space}
%  Set the catcode to 11 for use in private control sequence names.
\catcode`\@=11
%
%  Include all definitions related to the fonts msam, msbm and eufm, so that
%  when this file is used by itself, the results with respect to those fonts
%  are equivalent to what they would have been using AMS-TeX.
%  Most symbols in fonts msam and msbm are defined using \newsymbol;
%  however, a few symbols that replace composites defined in plain must be
%  defined with \mathchardef.

\def\undefine#1{\let#1\undefined}
\def\newsymbol#1#2#3#4#5{\let\next@\relax
 \ifnum#2=\@ne\let\next@\msafam@\else
 \ifnum#2=\tw@\let\next@\msbfam@\fi\fi
 \mathchardef#1="#3\next@#4#5}
\def\mathhexbox@#1#2#3{\relax
 \ifmmode\mathpalette{}{\m@th\mathchar"#1#2#3}%
 \else\leavevmode\hbox{$\m@th\mathchar"#1#2#3$}\fi}
\def\hexnumber@#1{\ifcase#1 0\or 1\or 2\or 3\or 4\or 5\or 6\or 7\or 8\or
 9\or A\or B\or C\or D\or E\or F\fi}

\font\tenmsa=msam10
\font\sevenmsa=msam7
\font\fivemsa=msam5
\newfam\msafam
\textfont\msafam=\tenmsa
\scriptfont\msafam=\sevenmsa
\scriptscriptfont\msafam=\fivemsa
\edef\msafam@{\hexnumber@\msafam}
\mathchardef\dabar@"0\msafam@39
\def\dashrightarrow{\mathrel{\dabar@\dabar@\mathchar"0\msafam@4B}}
\def\dashleftarrow{\mathrel{\mathchar"0\msafam@4C\dabar@\dabar@}}

\def\ulcorner{\delimiter"4\msafam@70\msafam@70 }
\def\urcorner{\delimiter"5\msafam@71\msafam@71 }
\def\llcorner{\delimiter"4\msafam@78\msafam@78 }
\def\lrcorner{\delimiter"5\msafam@79\msafam@79 }
\def\yen{{\mathhexbox@\msafam@55 }}
\def\checkmark{{\mathhexbox@\msafam@58 }}
\def\circledR{{\mathhexbox@\msafam@72 }}
\def\maltese{{\mathhexbox@\msafam@7A }}

\font\tenmsb=msbm10
\font\sevenmsb=msbm7
\font\fivemsb=msbm5
\newfam\msbfam
\textfont\msbfam=\tenmsb
\scriptfont\msbfam=\sevenmsb
\scriptscriptfont\msbfam=\fivemsb
\edef\msbfam@{\hexnumber@\msbfam}
\def\Bbb#1{{\fam\msbfam\relax#1}}
\def\widehat#1{\setbox\z@\hbox{$\m@th#1$}%
 \ifdim\wd\z@>\tw@ em\mathaccent"0\msbfam@5B{#1}%
 \else\mathaccent"0362{#1}\fi}
\def\widetilde#1{\setbox\z@\hbox{$\m@th#1$}%
 \ifdim\wd\z@>\tw@ em\mathaccent"0\msbfam@5D{#1}%
 \else\mathaccent"0365{#1}\fi}
\font\teneufm=eufm10
\font\seveneufm=eufm7
\font\fiveeufm=eufm5
\newfam\eufmfam
\textfont\eufmfam=\teneufm
\scriptfont\eufmfam=\seveneufm
\scriptscriptfont\eufmfam=\fiveeufm

%  Restore the catcode value for @ that was previously saved.
\csname amssym.def\endcsname

\expandafter\ifx\csname pre amssym.tex at\endcsname\relax \else \endinput\fi
%%     Otherwise we store the catcode of the @ in the csname.
\expandafter\chardef\csname pre amssym.tex at\endcsname=\the\catcode`\@
%%     Set the catcode to 11 for use in private control sequence
%%     names.
\catcode`\@=11
%  Most symbols in fonts msam and msbm are defined using \newsymbol.  A few
%  that are delimiters or otherwise require special treatment have already
%  been defined as soon as the fonts were loaded.  Finally, a few symbols
%  that replace composites defined in plain must be undefined first.
\newsymbol\boxdot 1200
\newsymbol\boxplus 1201
\newsymbol\boxtimes 1202
\newsymbol\square 1003
\newsymbol\blacksquare 1004
\newsymbol\centerdot 1205
\newsymbol\lozenge 1006
\newsymbol\blacklozenge 1007
\newsymbol\circlearrowright 1308
\newsymbol\circlearrowleft 1309
\undefine\rightleftharpoons
\newsymbol\rightleftharpoons 130A
\newsymbol\leftrightharpoons 130B
\newsymbol\boxminus 120C
\newsymbol\Vdash 130D
\newsymbol\Vvdash 130E
\newsymbol\vDash 130F
\newsymbol\twoheadrightarrow 1310
\newsymbol\twoheadleftarrow 1311
\newsymbol\leftleftarrows 1312
\newsymbol\rightrightarrows 1313
\newsymbol\upuparrows 1314
\newsymbol\downdownarrows 1315
\newsymbol\upharpoonright 1316
 \let\restriction\upharpoonright
\newsymbol\downharpoonright 1317
\newsymbol\upharpoonleft 1318
\newsymbol\downharpoonleft 1319
\newsymbol\rightarrowtail 131A
\newsymbol\leftarrowtail 131B
\newsymbol\leftrightarrows 131C
\newsymbol\rightleftarrows 131D
\newsymbol\Lsh 131E
\newsymbol\Rsh 131F
\newsymbol\rightsquigarrow 1320
\newsymbol\leftrightsquigarrow 1321
\newsymbol\looparrowleft 1322
\newsymbol\looparrowright 1323
\newsymbol\circeq 1324
\newsymbol\succsim 1325
\newsymbol\gtrsim 1326
\newsymbol\gtrapprox 1327
\newsymbol\multimap 1328
\newsymbol\therefore 1329
\newsymbol\because 132A
\newsymbol\doteqdot 132B
 
\newsymbol\triangleq 132C
\newsymbol\precsim 132D
\newsymbol\lesssim 132E
\newsymbol\lessapprox 132F
\newsymbol\eqslantless 1330
\newsymbol\eqslantgtr 1331
\newsymbol\curlyeqprec 1332
\newsymbol\curlyeqsucc 1333
\newsymbol\preccurlyeq 1334
\newsymbol\leqq 1335
\newsymbol\leqslant 1336
\newsymbol\lessgtr 1337
\newsymbol\backprime 1038
\newsymbol\risingdotseq 133A
\newsymbol\fallingdotseq 133B
\newsymbol\succcurlyeq 133C
\newsymbol\geqq 133D
\newsymbol\geqslant 133E
\newsymbol\gtrless 133F
\newsymbol\sqsubset 1340
\newsymbol\sqsupset 1341
\newsymbol\vartriangleright 1342
\newsymbol\vartriangleleft 1343
\newsymbol\trianglerighteq 1344
\newsymbol\trianglelefteq 1345
\newsymbol\bigstar 1046
\newsymbol\between 1347
\newsymbol\blacktriangledown 1048
\newsymbol\blacktriangleright 1349
\newsymbol\blacktriangleleft 134A
\newsymbol\vartriangle 134D
\newsymbol\blacktriangle 104E
\newsymbol\triangledown 104F
\newsymbol\eqcirc 1350
\newsymbol\lesseqgtr 1351
\newsymbol\gtreqless 1352
\newsymbol\lesseqqgtr 1353
\newsymbol\gtreqqless 1354
\newsymbol\Rrightarrow 1356
\newsymbol\Lleftarrow 1357
\newsymbol\veebar 1259
\newsymbol\barwedge 125A
\newsymbol\doublebarwedge 125B
\undefine\angle
\newsymbol\angle 105C
\newsymbol\measuredangle 105D
\newsymbol\sphericalangle 105E
\newsymbol\varpropto 135F
\newsymbol\smallsmile 1360
\newsymbol\smallfrown 1361
\newsymbol\Subset 1362
\newsymbol\Supset 1363
\newsymbol\Cup 1264
 
\newsymbol\Cap 1265
 
\newsymbol\curlywedge 1266
\newsymbol\curlyvee 1267
\newsymbol\leftthreetimes 1268
\newsymbol\rightthreetimes 1269
\newsymbol\subseteqq 136A
\newsymbol\supseteqq 136B
\newsymbol\bumpeq 136C
\newsymbol\Bumpeq 136D
\newsymbol\lll 136E
 
\newsymbol\ggg 136F
 
\newsymbol\circledS 1073
\newsymbol\pitchfork 1374
\newsymbol\dotplus 1275
\newsymbol\backsim 1376
\newsymbol\backsimeq 1377
\newsymbol\complement 107B
\newsymbol\intercal 127C
\newsymbol\circledcirc 127D
\newsymbol\circledast 127E
\newsymbol\circleddash 127F
\newsymbol\lvertneqq 2300
\newsymbol\gvertneqq 2301
\newsymbol\nleq 2302
\newsymbol\ngeq 2303
\newsymbol\nless 2304
\newsymbol\ngtr 2305
\newsymbol\nprec 2306
\newsymbol\nsucc 2307
\newsymbol\lneqq 2308
\newsymbol\gneqq 2309
\newsymbol\nleqslant 230A
\newsymbol\ngeqslant 230B
\newsymbol\lneq 230C
\newsymbol\gneq 230D
\newsymbol\npreceq 230E
\newsymbol\nsucceq 230F
\newsymbol\precnsim 2310
\newsymbol\succnsim 2311
\newsymbol\lnsim 2312
\newsymbol\gnsim 2313
\newsymbol\nleqq 2314
\newsymbol\ngeqq 2315
\newsymbol\precneqq 2316
\newsymbol\succneqq 2317
\newsymbol\precnapprox 2318
\newsymbol\succnapprox 2319
\newsymbol\lnapprox 231A
\newsymbol\gnapprox 231B
\newsymbol\nsim 231C
\newsymbol\ncong 231D
\newsymbol\diagup 231E
\newsymbol\diagdown 231F
\newsymbol\varsubsetneq 2320
\newsymbol\varsupsetneq 2321
\newsymbol\nsubseteqq 2322
\newsymbol\nsupseteqq 2323
\newsymbol\subsetneqq 2324
\newsymbol\supsetneqq 2325
\newsymbol\varsubsetneqq 2326
\newsymbol\varsupsetneqq 2327
\newsymbol\subsetneq 2328
\newsymbol\supsetneq 2329
\newsymbol\nsubseteq 232A
\newsymbol\nsupseteq 232B
\newsymbol\nparallel 232C
\newsymbol\nmid 232D
\newsymbol\nshortmid 232E
\newsymbol\nshortparallel 232F
\newsymbol\nvdash 2330
\newsymbol\nVdash 2331
\newsymbol\nvDash 2332
\newsymbol\nVDash 2333
\newsymbol\ntrianglerighteq 2334
\newsymbol\ntrianglelefteq 2335
\newsymbol\ntriangleleft 2336
\newsymbol\ntriangleright 2337
\newsymbol\nleftarrow 2338
\newsymbol\nrightarrow 2339
\newsymbol\nLeftarrow 233A
\newsymbol\nRightarrow 233B
\newsymbol\nLeftrightarrow 233C
\newsymbol\nleftrightarrow 233D
\newsymbol\divideontimes 223E
\newsymbol\varnothing 203F
\newsymbol\nexists 2040
\newsymbol\Finv 2060
\newsymbol\Game 2061
\newsymbol\mho 2066
\newsymbol\eth 2067
\newsymbol\eqsim 2368
\newsymbol\beth 2069
\newsymbol\gimel 206A
\newsymbol\daleth 206B
\newsymbol\lessdot 236C
\newsymbol\gtrdot 236D
\newsymbol\ltimes 226E
\newsymbol\rtimes 226F
\newsymbol\shortmid 2370
\newsymbol\shortparallel 2371
\newsymbol\smallsetminus 2272
\newsymbol\thicksim 2373
\newsymbol\thickapprox 2374
\newsymbol\approxeq 2375
\newsymbol\succapprox 2376
\newsymbol\precapprox 2377
\newsymbol\curvearrowleft 2378
\newsymbol\curvearrowright 2379
\newsymbol\digamma 207A
\newsymbol\varkappa 207B
\newsymbol\Bbbk 207C
\newsymbol\hslash 207D
\undefine\hbar
\newsymbol\hbar 207E
\newsymbol\backepsilon 237F
%  Restore the catcode value for @ that was previously saved.
\catcode`\@=\csname pre amssym.tex at\endcsname

%\endinput
%\input mssymb      % these must be input at the BEGINNING of this file; some
% MAXIMUM LINE LENGTH IS 80 CHARACTERS
%
%\font\eurofont=times at 10pt% contains Û symbol (as option-shift-2)
%
\font\fivebi=cmmib5
\font\fivebsy=cmbsy5
\font\sixrm=cmr6
\font\sixi=cmmi6
\font\sixbf=cmbx6
\font\sixsy=cmsy6
\font\sixmsa=msam5 at 6pt% This is needed since Bluesky's virtual fonts for
\font\sixmsb=msbm5 at 6pt% sizes 6.8.9 all come out too small (same below)
\font\sevenbi=cmmib7
\font\sevenbsy=cmbsy7
\font\eightrm=cmr8
\font\eightsl=cmsl8
\font\eightit=cmti8
\font\eighti=cmmi8
\font\eightbf=cmbx8
\font\eightsy=cmsy8
\font\eightmsa=msam7 at 8pt
\font\eightmsb=msbm7 at 8pt
\font\ninerm=cmr9
\font\ninesl=cmsl9
\font\nineit=cmti9
\font\ninei=cmmi9
\font\ninebf=cmbx9
\font\ninebi=cmmib10 scaled 900
\font\ninesy=cmsy9
\font\ninebsy=cmbsy10 scaled 900
\font\ninemsa=msam10 at 9pt
\font\ninemsb=msbm10 at 9pt
\font\tenbit=cmbxti10
\font\tenbsl=cmbxsl10
\font\tenbi=cmmib10
\font\tenbsy=cmbsy10
\font\twelvebf=cmbx12
\font\twelvebi=cmmib10 scaled 1200
\font\twelvebsy=cmbsy10 at 12pt

\let\sc=\sevenrm            % SMALL CAPS (in tenpoint)
\def\eightpoint{%
     %\font\eurofont=times at 8pt% contains Û symbol (as option-shift-2)
     \def\rm{\fam0\eightrm}%         see p.415
     \textfont0=\eightrm \scriptfont0=\sixrm \scriptscriptfont0=\fiverm
     \textfont1=\eighti \scriptfont1=\sixi \scriptscriptfont1=\fivei
     \textfont2=\eightsy \scriptfont2=\sixsy \scriptscriptfont2=\fivesy
     \textfont3=\tenex \scriptfont3=\tenex \scriptscriptfont3=\tenex
     \textfont\itfam=\eightit \def\it{\fam\itfam\eightit}%
     \textfont\slfam=\eightsl \def\sl{\fam\slfam\eightsl}%
     \textfont\bffam=\eightbf \scriptfont\bffam=\sixbf
     \scriptscriptfont\bffam=\fivebf \def\bf{\fam\bffam\eightbf}%
     \textfont\msbfam=\eightmsb \textfont\msafam=\eightmsa
     \scriptfont\msafam=\sixmsa \scriptfont\msbfam=\sixmsb
     \scriptscriptfont\msafam=\fivemsa \scriptscriptfont\msbfam=\fivemsb
      \skewchar\eighti='177 \skewchar\sixi='177
      \skewchar\eightsy='60 \skewchar\sixsy='60
     \normalbaselineskip=10pt%(normally 9pt; must be =height+depth in next line)
     \setbox\strutbox=\hbox{\vrule height7pt depth3pt width0pt}%(normally 7&2pt)
     \let\sc=\sixrm \let\big=\eightbig \normalbaselines\rm}
\def\ninepoint{%
     %\font\eurofont=times at 9pt% contains Û symbol (as option-shift-2)
     \def\rm{\fam0\ninerm}%         see p.415
     \textfont0=\ninerm \scriptfont0=\sixrm \scriptscriptfont0=\fiverm
     \textfont1=\ninei \scriptfont1=\sixi \scriptscriptfont1=\fivei
     \textfont2=\ninesy \scriptfont2=\sixsy \scriptscriptfont2=\fivesy
     \textfont3=\tenex \scriptfont3=\tenex \scriptscriptfont3=\tenex
     \textfont\itfam=\nineit \def\it{\fam\itfam\nineit}%
     \textfont\slfam=\ninesl \def\sl{\fam\slfam\ninesl}%
     \textfont\bffam=\ninebf \scriptfont\bffam=\sixbf
     \scriptscriptfont\bffam=\fivebf \def\bf{\fam\bffam\ninebf}%
     \textfont\msbfam=\ninemsb \textfont\msafam=\ninemsa
     \scriptfont\msafam=\sixmsa \scriptfont\msbfam=\sixmsb
     \scriptscriptfont\msafam=\fivemsa \scriptscriptfont\msbfam=\fivemsb
      \skewchar\ninei='177 \skewchar\sixi='177
      \skewchar\ninesy='60 \skewchar\sixsy='60
     \normalbaselineskip=11pt
     \setbox\strutbox=\hbox{\vrule height8pt depth3pt width0pt}%
     \let\sc=\sevenrm \let\big=\ninebig \normalbaselines\rm}
\def\tenpoint{%
     %\font\eurofont=times at 10pt% contains Û symbol (as option-shift-2)
     \def\rm{\fam0\tenrm}%         see p.415
     \textfont0=\tenrm \scriptfont0=\sevenrm \scriptscriptfont0=\fiverm
     \textfont1=\teni \scriptfont1=\seveni \scriptscriptfont1=\fivei
     \textfont2=\tensy \scriptfont2=\sevensy \scriptscriptfont2=\fivesy
     \textfont3=\tenex \scriptfont3=\tenex \scriptscriptfont3=\tenex
     \textfont\itfam=\tenit \def\it{\fam\itfam\tenit}%
     \textfont\slfam=\tensl \def\sl{\fam\slfam\tensl}%
     \textfont\bffam=\tenbf \scriptfont\bffam=\sevenbf
     \scriptscriptfont\bffam=\fivebf \def\bf{\fam\bffam\tenbf}%
     \textfont\msbfam=\tenmsb \textfont\msafam=\tenmsa
     \scriptfont\msafam=\sevenmsa \scriptfont\msbfam=\sevenmsb
     \scriptscriptfont\msafam=\fivemsa \scriptscriptfont\msbfam=\fivemsb
     \normalbaselineskip=12pt
%     \normalbaselineskip=13pt% (normally 12pt; here 13 = 9 + 4 of next line)
%     \setbox\strutbox=\hbox{\vrule height9pt depth4pt width0pt}%
     \let\sc=\sevenrm \let\big=\tenbig \normalbaselines\rm}
\catcode`@=11        % This allows the use of `@' in the next line (p.344)
\def\tenbig#1{{\hbox{$\left#1\vbox to8.5pt{}\right.\n@space$}}}
\def\ninebig#1{{\hbox{$\textfont0=\tenrm\textfont2=\tensy
      \left#1\vbox to7.25pt{}\right.\n@space$}}}
\def\eightbig#1{{\hbox{$\textfont0=\ninerm\textfont2=\ninesy
      \left#1\vbox to6.5pt{}\right.\n@space$}}}
\catcode`@=12        % This restores the `inhibiting' catcode of `@'.
\def\bold{%
     \textfont0=\tenbf \scriptfont0=\sevenbf \scriptscriptfont0=\fivebf
     \textfont1=\tenbi \scriptfont1=\sevenbi \scriptscriptfont1=\fivebi
     \textfont2=\tenbsy \scriptfont2=\sevenbsy \scriptscriptfont2=\fivebsy
       \textfont\itfam=\tenbit \def\it{\fam\itfam\tenbit}%
       \textfont\slfam=\tenbsl \def\sl{\fam\slfam\tenbsl}%
       \textfont\bffam=\tenbf \scriptfont\bffam=\sevenbf
       \textfont\msbfam=\tenmsb \textfont\msafam=\tenmsa
   \fam0\tenbf}
\def\bigbold{%
     \textfont0=\twelvebf \scriptfont0=\ninebf \scriptscriptfont0=\sevenbf
     \textfont1=\twelvebi \scriptfont1=\ninebi \scriptscriptfont1=\sevenbi
     \textfont2=\twelvebsy \scriptfont2=\ninebsy \scriptscriptfont2=\sevenbsy
     \fam0\twelvebf}
%
%
% MAXIMUM LINE LENGTH BELOW IS 80 CHARACTERS
%
% Instructions on how to use \item, and how to surround items by space,
% are given in paperformat.
%
\let\plainitem=\item% \item itself is redefined in paperformat.
\let\plainitemitem=\itemitem% \itemitem itself is redefined in paperformat.
%

%
% Û is option-shift-2 on the Mac
%
\def\openface{\Bbb}                %  needs mssymb
\def\N{{\openface N}}              %  in old files, \nat may be used instead

\def\Q{{\openface Q}}

\def\C{{\openface C}}
%
%
%                            SKIPS AND BREAKS
%
\def\g{\hskip.17em\relax}               %  breakable thin space
\def\th{\thinspace}                     %  non-breakable thin space
\def\nl{\hfil\break}
\newskip\Bigskipamount
   \Bigskipamount=2\baselineskip plus.5\baselineskip minus.3\baselineskip
\def\Bigbreak{\removelastskip\vskip0pt plus .1\vsize\penalty-1000
              \vskip0pt plus-.1\vsize\vskip\Bigskipamount}
% The following can be used after a display to remove the belowdisplayskip, e.g.
% when the display finishes a proof and is just followed by \endproof. If text
% follows, e.g. just a word like "and" leading to another display, one may wish
% to say \noskip\smallskip\noindent.
\def\noskip{\vskip-\lastskip\noindent}%use after display followed by short line
\def\Nobreak$$#1$${\postdisplaypenalty=10000$$#1$$\postdisplaypenalty=0}
%
%
%                             ABBREVIATIONS
%
     % \H will redefined below.

\let\doublebar=\| 
\def\|{\!\!\restriction\!\!}

  \let\sub=\sube

\def\supe{\supseteq}
  
\def\supne{\supsetneqq}
\def\sm{\smallsetminus}
\def\es{\emptyset}

 % italics, like K^n, P^n etc. Cf. I FOR "Inflated"

\def\:{\colon}
\def\minor{\preccurlyeq} 
\def\Minor{\succcurlyeq}

\def\slt{\mathrel{\hbox{$\minor$\kern-.6em\lower.33ex\hbox{${}_s\;$}}}}
\def\sgt{\mathrel{\mathchoice                        %("simplicial minors")
   {\hbox{$\Minor$\kern-.5em\lower.3ex\hbox{${}_s$}}}
   {\hbox{$\Minor$\kern-.5em\lower.3ex\hbox{${}_s$}}}
   {\hbox{$\scriptstyle\Minor\kern-.43em\lower.28ex\hbox{$\scriptstyle{}_s$}$}}
 {\hbox{$\scriptstyle\Minor\kern-.43em\lower.28ex\hbox{$\scriptstyle{}_s$}$}} }}    

\def\ucl(#1){\lfloor #1 \rfloor}% up-closure
\def\dcl(#1){\lceil #1 \rceil}% down-closure
%
%

%

% Example: '\interior P' puts a circle over the P.
% (the "\relax" is important: otherwise a capital letter A-F instead of P
% is interpreted as part of the number 7017,
% giving a "wrong math code complaint!!)
%
% The following commands can be used to 'specify' a relation (= #1) by putting
% something (= #2) underneath or above.
% Examples:\specrel<T, \specrel\sim G, \specrel={(1)}, \Specrel\Rightarrow?.
%
\def\specrel#1#2{\mathrel{\mathop{\kern0pt #1}\limits_{#2}}}
\def\Specrel#1#2{\mathrel{\mathop{\kern0pt #1}\limits^{#2}}}
%
% Next a version of \specrel for use in aligned equations;
% here the specification text is put in an \hbox of width 0,
% so as not to interfere with the alignment.
%
\def\alignspecrel#1#2{\mathrel{\mathop{\kern0pt #1}\limits_{\hbox
   to0pt{\hss$\scriptstyle#2$\hss}}}}
\def\alignSpecrel#1#2{\mathrel{\mathop{\kern0pt #1}\limits^{\hbox
   to0pt{\hss$\scriptstyle#2$\hss}}}}
\def\invlim{\specrel\lim{\raise 2pt\hbox{$\longleftarrow$}}}% inverse limit, projective limit
\def\proof{\removelastskip\penalty55\medskip\noindent{\bf Proof. }}
\def\noproof{{\unskip\nobreak\hfill\penalty50\hskip2em\hbox{}\nobreak\hfill%
       $\square$\parfillskip=0pt\finalhyphendemerits=0\par}\goodbreak}
\def\endproof{\noproof\bigskip}
% Syntax example for the following: \looseproof{Zweiter Beweis von Satz \xxxA}
%space follows in input

%
\newcount\refno
\def\ref#1#2\par{{\plainitem{[??]}#2\smallskip}}
\newtoks\thingtowrite %USED AGAIN LATER FOR INDEX
\long\def\writerefnumber#1{%
    \thingtowrite={#1}%
    \immediate\write\refnumbersfile{\the\thingtowrite}%
    }
\newwrite\refnumbersfile
\def\makerefnumbers{\immediate\openout\refnumbersfile=RefNumbers%
  \refno=0 \writerefnumber{\refno=0}
  \def\ref##1##2\par{\global\advance\refno by 1
    \writerefnumber{\global\advance\refno by 1 \newcounter##1 ##1=\the\refno}%
       % \newcounter replaces \newcount, which is forbidden inside a def.
       % When using the auto-generated file, say "\let\newcounter=\newcount"
       % before reading in that file.
    \plainitem{[\the\refno]}##2\smallskip% automatic numbering, ignoring ##1
    }%
  }
\def\autorefnumbers{\refno=0
  \def\ref##1##2\par{\advance\refno by 1\plainitem{[\the\refno]}##2\smallskip}
  }
\def\userefnumbers{\refno=0
  \def\ref##1##2\par{\advance\refno by 1\plainitem{[\the##1]}##2\smallskip}
  }
% TO USE, FIRST RUN WITH \makerefnumbers active (and no file named "RefNumbers" present
% in the directory to which such a file is to be written; it cannot overwrite).
% THEN DISABLE "\makerefnumbers", and make the following line active:
% \userefnumbers\let\newcounter=\newcount\input RefNumbers
%
%
\def\proclaimwithname #1. (#2) #3\par{{\bigbreak
  \clubpenalty=10000\noindent{\bf#1.\enspace}(#2)\nl
  {\sl #3}\par\bigbreak}}
\def\proposition (#1) #2\par{{\setbox0\hbox{(#1)\enspace}\bigbreak
   \sl\hangindent\the\wd0 \noindent\hskip\the\wd0
   \llap{\box0}\ignorespaces#2\par\bigbreak}}
\def\subsection #1\par{\vskip 3\medskipamount minus \smallskipamount\leftline{\bold #1}
        \penalty10000\smallskip\noindent}
      \def\section #1\par{\Bigbreak\centerline{\bf #1} % TO BE PHASED OUT
              \penalty10000\bigskip\noindent}
%
%         The following versions of the beginsection macro are for
%         section headings immediately followed by \proclaim:
\def\beginpsection #1\par{\Bigbreak\centerline{\bold #1}
        \penalty10000\bigskip}
\def\psubsection #1\par{\bigbreak\leftline{\bold #1}\penalty10000\bigskip}
%
%         The following macro positions its argument flush right,
%         like the \endproof box or an equation number.

%
%  The following three items usually need an \enditem. They are for use inside
%  \proclaim or another \item. Note that these items don't require their own
%  \par at the end, but they have no finishing \smallskip. (There seems to be
%  no way around this: we can't delimit by \par, since in the input this would
%  end the \proclaim.) So, if the item is followed by some text within the same
% \proclaim (say), the \enditem should (and can) be replaced with \smallskip.
%        
\def\pitem#1{\smallskip\advance\parindent by 3mm
             \plainitem{\rm(#1)}\advance\parindent by-3mm}
\def\pitemitem#1{\smallskip\advance\parindent by 3mm
             \plainitemitem{\rm(#1)}\advance\parindent by-3mm}
\def\varitemitem#1{{\setbox0\hbox{\hskip\parindent#1\enskip}
           \smallbreak\hangindent\the\wd0 \noindent\hskip\the\wd0
           \llap{#1\enskip}\ignorespaces}}
%
%        Note that the following variations of \item must end with \par.
%
\newdimen\newparindent
\def\iitem#1#2\par{\newparindent=\parindent \advance\newparindent by 3mm
           \smallbreak \hangindent\newparindent \noindent\hskip\newparindent
           \llap{{\rm #1}\enspace}\ignorespaces#2\par\smallbreak}
\def\iitemitem#1#2\par{\newparindent=\parindent \advance\newparindent by 3mm
           \smallbreak \hangindent2\newparindent \noindent\hskip2\newparindent
           \llap{{\rm #1}\enspace}\ignorespaces#2\par\smallbreak}
\def\varitem#1#2\par{{\setbox0\hbox{{\rm #1}\enspace}
           \smallbreak \hangindent\the\wd0 \noindent\hskip\the\wd0
           \llap{{\rm #1}\enspace}\ignorespaces#2\par\smallbreak}}
\def\enditem{\par}
% For use in Exercises and Hints
%
\def\Textindent#1{\par \advance\parindent by 3mm
                  \textindent{{\rm #1}} \advance\parindent by -3mm}
\def\indentedline#1{\advance\hsize by -\parindent \line{#1}
                   \advance\hsize by \parindent}
\def\iindentedline#1{\advance\parindent by 3mm
                     \advance\hsize by -\parindent
                     \line{#1}
                     \advance\hsize by \parindent
                     \advance\parindent by -3mm}
\newdimen\margin   % needed for macros \textdisplay & \ltextdisplay
%  The following macro takes 3 arguments, #1 and #3 in math-mode,
%  #2 as plain text. It displays #1\quad centered w.r.t. the whole \line
%  (if #2 leaves enough space), adds \quad#2 to the right of #1\quad, and puts
%  #3 flush right. The arguments must be separated by &'s. The argu-
%  ments themselves may be empty, but there must be 2 ampersands.
%  Examples: $$\textdisplay x=y &for all $x\in X$& (1')$$
%            $$\textdisplay x=y &for all $x\in X$&$$
\def\textdisplay#1&#2&#3$${\margin=\hsize
          \setbox1=\hbox{$\displaystyle#1\quad$}%
          \setbox2=\hbox{\quad#2\qquad$#3$}%
                     \advance\margin by-\wd1
                     \divide\margin by 2
   \ifdim\wd2 < \margin
      \box1\rlap{\quad#2}\eqno#3$$%
   \else
      \line{\qquad\hfil \box1\quad #2 \qquad $#3$}$$%
   \fi}
%
%  The following macro is the `\leqno' version of \textdisplay; argument
%  #1 is the \leqno, #2 is the formula to be displayed, and #3 is the
%  text following the formula.
\def\ltextdisplay#1&#2&#3$${\margin=\hsize
           \setbox2=\hbox{$\displaystyle#2\quad$}
           \setbox3=\hbox{\quad#3\qquad}
                     \advance\margin by-\wd2
                     \divide\margin by 2
   \ifdim\wd3 < \margin
      \line{$#1$\hfil\box2\hbox to \margin{\box3\hfil}}$$%
   \else
      \line{$#1$\qquad\hfil\box2\quad #3\qquad} $$%
   \fi}
%
%   The next macro displays and centres #1, typically a paragraph of text.
%   #2, typically an `eqno', is set flush right, vertically centred, and
%   has to be in math mode.
%   Example: \textno This is a ... lot of text&(3')\par
\def\textno#1&#2\par{%
   \margin=\hsize
   \advance\margin by -4\parindent
          \setbox1=\hbox{\sl#1}%
   \ifdim\wd1 < \margin
      $$\box1\eqno#2$$\endgraf%
   \else
      \bigbreak
      \line{\indent$\vcenter{\advance\hsize by -3\parindent
      \sl\noindent#1}\hfil#2$}%
      \bigbreak
   \fi}
%
%   The same with a left eqno:
\def\textlno#1&#2\par{%
   \margin=\hsize
   \advance\margin by -4\parindent
          \setbox1=\hbox{\sl#1}%
   \ifdim\wd1 < \margin
      $$\box1\leqno#2$$%
   \else
      \bigbreak
      \line{$#2\hfil\vcenter{\advance\hsize by -3\parindent
          \sl\noindent#1}\hskip\parindent$}%
      \bigbreak
   \fi}
%
%
%                        MARGINAL HACKS ETC.
%
\newcount\commentno
\def\COMMENT#1{$^{<\the\commentno>}$%
     \vadjust{\vbox to 0pt{\vss\vskip-8pt\rightline{%
     \rlap{\hbox{\hskip7mm \vbox{\pretolerance=-1
     \doublehyphendemerits=0 \finalhyphendemerits=0
     \hsize40mm\tolerance=10000\eightpoint
     \lineskip=0pt\lineskiplimit=0pt
     \rightskip=0pt plus16mm\baselineskip8pt\noindent
     \hskip0pt       %(without this, the first word is never hyphenated!)
     {$\langle$\the\commentno. #1$\rangle$}\endgraf}}}}\vss}}%
     \global\advance\commentno by1}%
\def\writecommentsasfootnotes{%
 \def\COMMENT{\global\advance\commentno by1\footnote{$^{<\the\commentno>}$}}%
 }
\def\nocomments{\def\COMMENT##1{}}
%
%
%   The \? macro puts
%   the argument #1 in the left margin. Examples: \??, \?{What nonsense!}.
\def\?#1{\vadjust{\vbox to 0pt{\vss\vskip-8pt\leftline{%
     \llap{\hbox{\vbox{\pretolerance=-1
     \doublehyphendemerits=0\finalhyphendemerits=0
     \hsize16truemm\tolerance=10000\eightpoint
     \lineskip=0pt\lineskiplimit=0pt
     \rightskip=0pt plus16truemm\baselineskip8pt\noindent
     \hskip0pt        %(without this, the first word is never hyphenated!)
     #1\endgraf}\hskip7truemm}}}\vss}}}
\def\d{}% DISABLES WHAT'S JUST BEEN DEFINED ABOVE
%
% The following \ds macro is the "silent" version of \d: it writes its argument
% in the margin (or index file), just as \d does, but not into the current text.
% Note the different syntax: no delimiting blank in input.
%
%  CHANGE: I got fed up with marginal reminders to define things, so this
%  feature is disabled right away by the line "\def\ds#1{}" below. Note that
%  \makeindex still works, because it defines \d and \ds anew.
%
%\def\ds#1{\ifmmode
%     \vadjust{\vbox to 0pt{\vss\vskip-8pt\leftline{%
%     \llap{\hbox{\vbox{\pretolerance=-1
%    \doublehyphendemerits=0\finalhyphendemerits=0
%     \hsize16truemm\tolerance=10000\eightpoint
%     \lineskip=0pt\lineskiplimit=0pt
%     \rightskip=0pt plus16truemm\baselineskip8pt\noindent
%     define $#1$!\endgraf}\hskip7truemm}}}\vss}}%
%   \else
%     \vadjust{\vbox to 0pt{\vss\vskip-8pt\leftline{%
%     \llap{\hbox{\vbox{\pretolerance=-1
%     \doublehyphendemerits=0\finalhyphendemerits=0
%     \hsize16truemm\tolerance=10000\eightpoint
%     \lineskip=0pt\lineskiplimit=0pt
%     \rightskip=0pt plus16truemm\baselineskip8pt\noindent
%     \hskip0pt        %(without this, the first word is never hyphenated!)
%     define ``#1''!\endgraf}\hskip7truemm}}}\vss}}%
%   \fi}
%
\def\ds#1{}% DISABLES WHAT'S JUST BEEN DEFINED ABOVE
%
% The following is a device (due to CET1) which allows to \write something to
% a file without expanding all the tokens completely. (This would result in
% the use of lots of "@"s, which make the file untexable. Another way around the
% problem would be to precede every argument of \write by "\catcode`@=11" (and
% to set it back to 12 after the argument). That would make the file texable,
% but since it would still be difficult to read, the method below is better.
% NOTE: This works well for a Remarks file (say), but NOT for an index file. 
% The reason is that in order to get the page numbers of entries from the first
% paragraph of a page right, one has to say "\write" rather than
% "\immediate\write". But in the def of \indexwrite one has to say "\immediate",
% since otherwise the token gets overwritten and rather
% than n index words from a page the index will
% contain the last index word of that page n times. But if the index words are
% written to file immediately while their page numbers get queued, then the two
% get separated from one another.
%
%\newtoks\thingtowrite % NOW EARLIER
\long\def\indexwrite#1{%
    \thingtowrite={#1}%
    \immediate\write\index{\the\thingtowrite}%
    }
%
%  The following macro \makeindex redefines \d and \ds.
%  It suppresses the appearance of \d's arguments in the margin
%  and writes them to a file called "index" instead.
%
\newwrite\index
\def\makeindex{\immediate\openout\index=index%
   \immediate\write\index{\catcode`@=11}%
   \def\d##1 {\ifmmode
     \write\index{$##1$, }%
     \write\index{\the\count0}\write\index{}% blank line for para
   \else
     \write\index{{##1}, }%
     \write\index{\the\count0}\write\index{}% blank line for para
   \fi {##1} }% It's important to put the text AFTER the def: otherwise
              % the blank is not merged with possible blanks following in
              % the input file, which may result in an additional blank line
      \def\ds##1{\ifmmode
     \write\index{$##1$, }%
     \write\index{\the\count0}\write\index{}% blank line for para
   \else
     \write\index{##1, }%
     \write\index{\the\count0}\write\index{}% blank line for para
   \fi}}
%
%   The \m macro puts its (direct!) argument in the right margin,
%   leaving it also in the text (in italics). 
%   Syntax example: this \m{word\/} is defined here
%   and will be shown in the margin; similarly the set $X$ in $A =: \m X$.
%   Intended use: in preprints, at places where something is first defined
%   (to enhance readability of the paper, in the absence of an index).
%   To disable, say "\disablems" at beginning of file.
\newdimen\gap% gap between text and hack
\gap=3truemm
\newdimen\hackwidth
\hackwidth=15truemm
% NOT simply "{}": when \m is active
 % it puts a strut in the line, so this should be here also when \m is disabled
 % to avoid a change in vertical space (and hence possibly in page breaks)
% \ds writes to index file, \mo in margin and text
\def\mo#1{\ifmmode {#1}\else {\it#1}\fi\mos{#1}}
% \ds writes to index file, \mos in margin
\def\mos#1{\ifmmode
     \strut\vadjust{\vbox to 0pt{\vss\kern-11pt\leftline{%
     \llap{\hbox{\vbox{\pretolerance=-1
     \doublehyphendemerits=0\finalhyphendemerits=0
     \hsize\hackwidth\tolerance=10000\eightpoint
     \lineskip=0pt\lineskiplimit=0pt
     \rightskip=0pt plus\hsize\baselineskip8pt\noindent
     $#1$\strut\endgraf}\hskip\gap }}}\vss}}%
   \else
     \strut\vadjust{\vbox to 0pt{\vss\kern-11pt\leftline{%
     \llap{\hbox{\vbox{\pretolerance=-1
     \doublehyphendemerits=0\finalhyphendemerits=0
     \hsize\hackwidth\tolerance=10000\eightpoint
     \lineskip=0pt\lineskiplimit=0pt
     \rightskip=0pt plus\hsize\baselineskip8pt\noindent
     \hskip0pt    %(without this, the first word is never hyphenated!)
     {\sl#1}\strut\endgraf}\hskip\gap }}}\vss}}%
   \fi}%
\newcount\remarkno
\def\REMARK#1{{\footnote{${}^{\the\remarkno}$}{{#1}}%
   \global\advance\remarkno by1}}

%
%
%
%                              PICTURES
%
\def\picture #1 by #2 (#3){
  \vbox to #2{
          \vfill
          \special{picture #3}
          \hrule width #1 height 0pt depth 0pt
           }}
\newdimen\topfiguremargin
   \topfiguremargin=0pt                                  % default
\newdimen\bottomfiguremargin
   \bottomfiguremargin=\medskipamount                    % default
\newdimen\normalpictureheight
\normalpictureheight=40mm
%   The following macro uses TeXtures pictures; these MUST be named
%   Fig.1, Fig.2.5 etc. (without a space after the '.'), as in the macro call
%itself. The width (#2) and height (#3) should have the original values
%of the TeXtures picture, to facilitate proper scaling. The heightfactor (#4)
%is divided by 1000 and used to scale the \normalpictureheight. Thus, height-
%factor 500 (2000) results in a picture of half (twice) the normalpicureheight.
%   Syntax example for a figure at the standard \normalpictureheight:
%\Fig.2 (538pt by 536pt; heightfactor: 1000; caption: This is a short caption)
\def\Fig.#1 (#2by#3; heightfactor:#4; caption:#5) {{%
   \dimen2=\normalpictureheight
   \dimen0=#2                          % computing width
      \divide\dimen2 by 1000
      \multiply\dimen2 by#4              % \dimen2 := intended pictureheight
   \count2=\dimen2                  % computing scalefactor
      \dimen1=#3                             % \dimen1 := actual pictureheight
   \count1=\dimen1
   \divide\count1 by 1000
   \divide\count2 by \count1          % \count2 := scalefactor (times 1000)
%      \message{scalefactor in Fig.#1 is \the\count2}%
   \divide\dimen0 by 1000
   \multiply\dimen0 by \count2      % \dimen0 := width
         \dimen1=\hsize
         \advance\dimen1 by -\dimen0
         \divide\dimen1 by 2               % \dimen1 := margin
   \midinsert
   \vbox to \topfiguremargin{\vfil}
   \noindent\hskip\dimen1
   \picture\dimen0 by \dimen2  (Fig.#1 scaled \the\count2)%
   \vskip\bottomfiguremargin                     % beginning caption
      \ninepoint
      \parindent=.1\hsize\narrower\narrower
      \setbox0\hbox{#5}
      \ifdim\wd0 < .6\hsize
           \centerline{F{\sc IGURE} #1.\hskip1em#5}
       \else
           \plainitem{F{\sc IGURE} #1. }#5\par
       \fi
   \vskip0pt\endinsert}}
%
%The following \textpicture macro is for inserting pictures in the current
%text line. The width (#2) and height (#3) should have the original values of
%the TeXtures picture, to facilitate proper scaling. #4 offers the opportunity
%to reserve vertical space for the picture by saying "height15pt depth10pt" or
%so; this will be the height of the \vbox containing the picture (default=0).
%#6 is the amount by which the bottom left corner is placed below the baseline
%(poss.neg.), #5 is the horizontal extension of the picture. Syntax examples:
%   \textpicture flower(538pt by 536pt; width5em lower-5pt)
%   \textpicture bug(538pt by 536pt; height20pt depth20pt width5em lower20pt)
\def\textpicture #1(#2by#3; #4width#5lower#6){{%
  % computing scalefactor
      \dimen0=#5\count2=\dimen0                    % desired width
      \dimen0=#2\count1=\dimen0                    % actual width
   \divide\count1 by 1000
   \divide\count2 by \count1                 % \count2 is now = scalefactor
  %\count3=11820\divide\count3 by \count2
  %\message{The vertices of #1 should have width \the\count3}
   \hbox{\vrule #4width0pt\vbox to 0pt{\vss\vskip#6%
      \special{picture #1 scaled \the\count2}\hrule width#5 height0pt\vss}}}}
%
%
%The following \figure macro builds on the \epsf macro and includes an EPS file.
% THE FOLLOWING LINE HAS TO OCCUR AT THE START OF THE MAIN FILE TO BE TEXED:
% \input epsf.def
%
% #1 is just a figure number to be used in the caption.
% #2 is the caption itself.
% #3 is the file name of the figure; this must an EPS file
%    (or PS with bounding box).
% #4 is for scaling, but disabled now (see below)
%
% Syntax example: \figure 3. The funny graph $G$ (Funny.graph.eps; 800)
%
\def\figure #1. #2 (#3; #4) {{%
   \def\bigskip{\par\ifdim\lastskip<\bigskipamount\removelastskip
                                              % eg. for abb after \endproof
      \vskip\bigskipamount\fi}% takes effect inside the def(!) of midinsert
   \midinsert\vskip\topfiguremargin
   \dimen0=\normalpictureheight
      \divide\dimen0 by 1000
      \multiply\dimen0 by#4        % \dimen0 := intended pictureheight
   \centerline{\epsfbox{#3.eps}}%                  good placing - use this!
   \vskip\bottomfiguremargin                     % beginning caption
      \ninepoint
      \parindent=.1\hsize\narrower\narrower
      \setbox0\hbox{#2}
      \ifdim\wd0 < .6\hsize
           \centerline{F{\sc IGURE} #1.\hskip1em#2}
       \else
           \plainitem{F{\sc IGURE} #1. }#2\par
       \fi
  \endinsert}}
%
% The following is for notes for talks to be hidden on transparancies.
% Say \hide{bla} to make "bla" disappear (but it will take up space);
% \showhidden will show all hidded text in grey (eg for paper notes).
%

%
%                               REFERENCES
%
\def\Abh#1 {{\sl Abh.\g Math.\g Sem.\g Univ.\g Hamburg\penalty100\ \bf#1\ }}
\def\AMASH#1 {{\sl Acta Math.\g Acad.\g Sci.\g Hung.\penalty100\ \bf#1\ }}
\def\Advances#1 {{\sl Adv.\g Math.\penalty100\ \bf#1\ }}
\def\Annals#1 {{\sl Ann.\g Math.\penalty100\ \bf#1\ }}
\def\AnnComb#1 {{\sl Ann.\g Comb.\penalty100\ \bf#1\ }}
\def\AMM#1 {{\sl Amer.\g Math.\g Monthly\penalty100\ \bf#1\ }}
\def\Archiv#1 {{\sl Arch.}\g {\sl Math.\penalty100\ \bf#1\ }}
\def\ArsComb#1 {{\sl Ars Comb.\penalty100\ \bf#1\ }}% FROMERLY \AC
\def\CJM#1 {{\sl Can.\g J.\th Math.\penalty100\ \bf#1\ }}
\def\Comb#1 {{\sl Com\-bi\-na\-to\-ri\-ca\penalty100\ \bf#1\ }}
\def\CPC#1 {{\sl Comb.\g Probab.\g Comput.\penalty100\ \bf#1\ }}
\def\Crelle#1 {{\sl J.}\th {\sl Reine Angew.}\g
    {\sl Math.\penalty100\ \bf#1\ }}
\def\DM#1 {{\sl Discrete Math.\penalty100\ \bf#1\ }}
\def\DAM#1 {{\sl Discrete Appl.\g Math.\penalty100\ \bf#1\ }}
\def\EJC#1 {{\sl Eur.}\g{\sl J.}\g{\sl Comb.\penalty100\ \bf#1\ }}
\def\EJ#1 {{\sl Electronic.}\g{\sl J.}\g{\sl Comb.\penalty100\ \bf#1\ }}
\def\GC#1 {{\sl Graphs Comb.\penalty100\ \bf#1\ }}
\def\IJ#1 {{\sl Isr.\g J.\th Math.\penalty100\ \bf#1\ }}
\def\Inv#1 {{\sl In\-vent.\g math.\penalty100\ \bf#1\ }}
\def\JAlg#1 {{\sl J.}\th {\sl Algorithms\penalty100\ \bf#1\ }}
\def\JCTA#1 {{\sl J.}\th {\sl Comb.}\g {\sl Theory~A\penalty100\ \bf#1\ }}
\def\JCTB#1 {{\sl J.}\th {\sl Comb.}\g {\sl Theory~B\penalty100\ \bf#1\ }}
\def\JGT#1 {{\sl J.}\th {\sl Graph Theory\penalty100\ \bf#1\ }}
\def\BLMS#1 {{\sl Bull.\g Lond.\g Math.\g Soc.\penalty100\ \bf#1\ }}
\def\JLMS#1 {{\sl J.\g Lond.\g Math.\g Soc.\penalty100\ \bf#1\ }}
\def\PLMS#1 {{\sl Proc.\g Lond.\g Math.\g Soc.\penalty100\ \bf#1\ }}
\def\Order#1 {{\sl Order\ \bf#1\ }}
\def\Random#1 {{\sl Random Struct.\g Alg.\penalty100\ \bf#1\ }}
\def\MA#1 {{\sl Math.}\g {\sl Ann.\penalty100\ \bf#1\ }}
\def\MN#1 {{\sl Math.}\g {\sl Nachr.\penalty100\ \bf#1\ }}
\def\MPCPS#1 {{\sl Math.\g Proc.\g Camb.\g Phil.\g Soc.\penalty100\ \bf#1\ }}
\def\MS#1 {{\sl Math.}\g {\sl Scand.\penalty100\ \bf#1\ }}
\def\MZ#1 {{\sl Math.}\g {\sl Zeit.\penalty100\ \bf#1\ }}
\def\BAMS#1 {{\sl Bull.\th Amer.\g Math.\g Soc.\penalty100\ \bf#1\ }}
\def\JAMS#1 {{\sl J.\th Amer.\g Math.\g Soc.\penalty100\ \bf#1\ }}
\def\MAMS#1 {{\sl Mem.\g Amer.\g Math.\g Soc.\penalty100\ \bf#1\ }}
\def\PAMS#1 {{\sl Proc.\g Amer.\g Math.\g Soc.\penalty100\ \bf#1\ }}
\def\SIAM#1 {{\sl SIAM J.\g Discrete Math.\penalty100\ \bf#1\ }}
\def\SLNM#1 {{\sl Springer Lecture Notes in Mathematics\penalty100\ \bf#1\ }}
\def\TAMS#1 {{\sl Trans.\g Amer.\g Math.\g Soc.\penalty100\ \bf#1\ }}
\def\TCSA#1 {{\sl Theor.\g Comput.\g Sci.~A\penalty100\ \bf#1\ }}
%
%
%
%
%
%                     ONLY FOR TeXtures without mssymb.tex:
%
%\def\N{{\rm \rlap I{\kern.18em}N}}   
%  \catcode`@=11        % This allows the use of `@' in the next line (p.344)
%\def\not#1{\mathrel{\mathpalette\c@ncel#1}} % use (only) with Imagewriter.
%  \catcode`@=12        % This restores the `inhibiting' catcode of `@'.
%\def\subsetneqq{\mathrel %{\mathchoice
%   {\vcenter{
%      \hbox{\lower6pt\hbox{$\scriptstyle\subset$}}
%      \hbox{\raise3pt\hbox{$\flatneq$}}}} }
%   \def\flatneq{\rlap {$\scriptstyle =$} {\kern1.5pt} {\scriptscriptstyle /}}
%\def\square{\Square53}            % See def. of \Square below.
%\def\Square#1#2{{\vbox{\hrule height.#2pt
%       \hbox{\vrule width.#2pt height#1pt \kern#1pt
%          \vrule width.#2pt}
%       \hrule height.#2pt}}}
%\def\nexists{\hbox{\rm\raise1.2pt\rlap/$\exists$}}
%
%
%
%
%
%
%\def\language=#1{}% For use with Textures versions < 1.3
%
%
%                  MODIFICATIONS TO PLAIN TeX
%
%\catcode`[=\active \catcode`]=\active
%  \def[{\thinspace\lbrack\thinspace}
%  \def]{\thinspace\rbrack}
%\def\{{\lbrace\thinspace}
%\def\}{\thinspace\rbrace}
%
\bigskipamount=1\baselineskip plus.3\baselineskip minus.3\baselineskip
\medskipamount=\bigskipamount\divide\medskipamount by 2
\smallskipamount=\medskipamount\divide\smallskipamount by 2 % (p.349)
\medmuskip = 3mu plus 2mu minus 1mu
\thickmuskip = 6mu plus 4mu minus 2mu % (for spacing in formulae; pp.168/170)
\def\smallbreak{\par \ifdim\lastskip<\smallskipamount
   \removelastskip \penalty-100 \smallskip \fi}
\def\medbreak{\par \ifdim\lastskip<\medskipamount
   \removelastskip \penalty-250 \medskip \fi}
\def\bigbreak{\par \ifdim\lastskip<\bigskipamount
   \removelastskip \penalty-500 \bigskip \fi}
\catcode`@=11        % This allows the use of `@' in the next line (p.344)
  \def\raggedbottom{\topskip10pt plus20pt \r@ggedbottomtrue} % The amount of
%                      permitted raggedness is controlled by the `plus' item.
\catcode`@=12        % This restores the `inhibiting' catcode of `@'.
\def\ge{\geqslant}% \geq remains available for the default version of \ge.
\def\le{\leqslant}% \leq remains available for the default version of \le.
\let\elt=\in
\def\in{\mathrel{\mathchoice
   {\raise .7pt \hbox{$\scriptstyle\elt$}}
   {\raise .7pt \hbox{$\scriptstyle\elt$}}
   {\raise .5pt \hbox{$\hskip .5pt\scriptscriptstyle\elt\hskip .5pt$}}
   {\raise.35pt \hbox{$\scriptscriptstyle\elt$}} }}
\let\hasaselt=\owns
\def\owns{\mathrel{\mathchoice
   {\raise .7pt \hbox{$\scriptstyle\hasaselt$}}
   {\raise .7pt \hbox{$\scriptstyle\hasaselt$}}
   {\raise .5pt \hbox{$\hskip .5pt\scriptscriptstyle\hasaselt\hskip .5pt$}}
   {\raise.35pt \hbox{$\scriptscriptstyle\hasaselt$}} }}
\let\exis=\exists
   \def\exists{\exis\>}
\let\nexis=\nexists
   \def\nexists{\nexis\>}
                            % To be phased out
\let\foral=\forall
   \def\forall{\foral\>}
\let\Rightarro=\Rightarrow
   \def\Rightarrow{\>\Rightarro\>}
\let\mi=\min
   \def\min{\mi\>}
\let\ma=\max
   \def\max{\ma\>}
\let\su=\sup
   \def\sup{\su\>}
\let\inff=\inf
   \def\inf{\inff\>}
\mathchardef\to="2221   % = \rightarrow, but of `binop' type (p.154)
\def\proclaim #1.#2 #3\par{\bigbreak
   \noindent{\bf#1.}#2\enspace{\sl#3}\par\bigbreak}
% The second argument above is optional and intended for references. Note that
% it is terminated by a space, which must therefore be present unsuppressed in
% input. Syntax examples: "\proclaim Thm {1.3}. This is the theorem.\par" or
% "\proclaim Thm \Euler.\five{} This is Euler's theorem.\par", where \Euler
%  expands to {1.3} and \five to \th [5], say. (Note that \proclaim Thm 1.3.
% would treat the 3. as argument #2 (incorrectly) and fail to set it in bold.)
% Or directly: "\proclaim Theorem \xxxVTop.~[\the\ref] blabla" (note the ~).
%
\newskip\sectionheadlineskipamount
\sectionheadlineskipamount=8pt plus 2pt minus 1pt
\def\beginsection #1\par{\Bigbreak\centerline{\bold #1}
        \penalty10000\vskip\sectionheadlineskipamount\noindent}
\let\ffootnote=\footnote
\def\footnote#1#2{\ffootnote{#1}{\eightpoint#2\vskip-12pt}}
%                  (The \vskip inserts an implicit \par, which has two
%                    effects: first, the desired effect of wrapping
%                    up the last paragraph of the footnote giving it the
%                    correct linespacing (that of \eightpoint), secondly
%                    the undesired effect of starting a new paragraph with
%                    a strut. To counteract the arising blank vertical
%                    space, the skip is chosen negative.)
%
% Note that the following redefinition of \item finishes with a smallbreak.
% Since two smallbreaks in sequence result in only a single smallbreak, this
% gives a smallskip at the beginning, between any two items, and at the end.
% if additional space is desired before and after a series of items, say
% \medskip before and \par\smallskip (not \smallbreak) after the series of
% items.
%
\def\item#1#2\par{\parindent=10mm\smallbreak\hang\indent
                  \llap{{\rm #1}\enspace}\ignorespaces#2\par\smallbreak
                  \parindent=7mm}
\def\itemitem#1#2\par{\parindent=10mm\smallbreak
                  \indent\hangindent2\parindent\indent
                  \llap{{\rm #1}\enspace}\ignorespaces#2\par\smallbreak
                  \parindent=7mm}
%
%
%                     FORMAT PARAMETERS
%
\pretolerance=0 %This prevents line breaks in maths formulas - I don't know why.
\tolerance=2000
%\fontdimen2\tenrm=3.8pt
%\fontdimen2\tensl=3.8pt
%\fontdimen2\tenit=4pt
\baselineskip=13pt                 %(DEFAULT IS 12pt)
\vsize=200mm                   % (used to be 240truemm; changed 06/98)
\hsize=120mm                   % (used to be 140truemm; changed 06/98)
\hoffset=9mm                   % (used to be 9truemm; changed 06/98)
\parindent=7mm
\relpenalty=2000 \binoppenalty=5000  % DISCOURAGES BREAKS IN FORMULAS
\hyphenpenalty=100
\abovedisplayskip=12pt plus3pt minus4pt
\belowdisplayskip=12pt plus3pt minus4pt    % (p.348)
%
%   Unfortunately, TeX is unable to adjust the skip following a display to the
%   length of the line **below** the display, in the way in which it chooses
%  between \abovedisplayskip and \abovedisplayshortskip depending on the length
%   of the line above it. Thus, such adjustment has to be done by hand: setting
%
\belowdisplayshortskip=9pt plus3pt minus3pt
    % (formerly 12pt, like \belowdisplayskip)
%
%   reduces the effect of the standard 'short' version, while
%
%                  \def\noskip{\vskip-\lastskip\noindent}
%
%   (which is contained in macros.tex) removes any belowdisplay skip
%   (long or short) altogether.
%   One may want to say \smallskip\noindent just after it.
%
   %\language=\Germanlanguage %Commented out to work with (non-bilingual) Plain, for ArXiv
% DON'T SAY "\German" HERE; THE READING-IN OF GERMAN.TEX WHICH THIS
% CAUSES CREATES ERRORS E.G. WITH LATER USE OF THE " CHARACTER
% IN THE \specialS USED FOR MAKING PDF LINKS, BECAUSE " BECOMES
% ACTIVE AND IS DEFINED IN AN UNDESIRED WAY.
% ***HOWEVER: if, for some reason, it becomes necessary to say \German
% here, it is possible to toggle the "-character between active and normal,
% using the commands \mdqon and \mdqoff defined in German.tex.
%
 % UNFORTUNATELY, NO UMLAUTS CAN BE USED IN \HYPHENATION 
 \hyphenation{Baum-ord-nung Baum-ord-nun-gen End-ecke End-ecken kur-zen
Kur-zen Graphen-ei-gen-schaft Graphen-ei-gen-schaften he-raus he-raus-ar-bei-ten
he-raus-zu-ar-bei-ten Schnitt-raum}%
 %\English %Commented out to work with (non-bilingual) Plain, for ArXiv
 \hyphenation{ac-cess-ible ana-log-ous ana-log-ous-ly ana-lyze ana-lyse
ana-ly-sis answer answers aver-age bundle bundles Buch-ge-sell-schaft col-our
col-ours col-oured col-our-ing col-our-ings con-struct-ible con-struct-ive
con-struct-ive-ly co-rol-lary Co-rol-lary des-cend des-cend-ing Deut-sche
end-li-cher de-fi-ni-tion de-fi-ni-tions De-fi-ni-tion equi-val-ent
equi-val-ence Euler-ian exist-ence every Gra-phen Hamil-ton-ian homeo-mor-phic
homeo-mor-phism homeo-mor-phisms hy-po-thesis hy-po-theses in-ac-cess-ible
ir-regu-lar ir-regu-lar-ity method methods modi-fi-ca-tion mono-chro-matic par-ticu-lar
pro-po-si-tion pro-po-si-tions Pro-po-si-tion regu-lar regu-lar-ity regu-lar-ly
sig-ni-fi-cant sig-ni-fi-cant-ly sig-ni-fi-cance to-po-lo-gical to-po-lo-gical-ly
un-at-tached un-end-li-cher using Using Wis-sen-schaft-li-che}
\userefnumbers\let\newcounter=\newcount\refno =0
\global \advance \refno by 1 \newcounter \refSpanningTrees  \refSpanningTrees =\the \refno 
\global \advance \refno by 1 \newcounter \refBook  \refBook =\the \refno 
\global \advance \refno by 1 \newcounter \refTopSurvey  \refTopSurvey =\the \refno 
\global \advance \refno by 1 \newcounter \refMayaEndDeg  \refMayaEndDeg =\the \refno 
\global \advance \refno by 1 \newcounter \refMayaBanff  \refMayaBanff =\the \refno 
\global \advance \refno by 1 \newcounter \refMayaZamora  \refMayaZamora =\the \refno 
\global \advance \refno by 1 \newcounter \refHalinEnd  \refHalinEnd =\the \refno

\hbox{}\vskip2pt
\centerline{\bigbold Forcing finite minors in sparse infinite graphs}\smallskip
\centerline{\bigbold by large-degree assumptions}\vskip 4mm
\centerline{Reinhard Diestel}

\advance\vsize by 1cm
\nocomments

\def\xxxNaiveThm{1}% Theorem
\def\xxxCtble{2}% Thm
\def\xxxUnctble{3}% Thm
\def\xxxPos{4}% Thm
\def\xxxCtbleTop{5}% Corollary

\def\C{{\cal C}}

\def\beginsection #1\par{\vskip2\bigskipamount\bigbreak\centerline{\bold #1}
        \penalty10000\vskip\sectionheadlineskipamount\noindent}

\bigskip\medskip
{\narrower\narrower\ninepoint\noindent
   Developing further Stein's recent notion of relative end degrees in infinite graphs, we investigate which degree assumptions can force a locally finite graph to contain a given finite minor, or a finite subgraph of given minimum degree. This is part of a wider project which seeks to develop an extremal theory of sparse infinite graphs. \par}

%\vskip-\medskipamount\vskip0pt

\beginsection Introduction

Perhaps the most basic question in extremal graph theory asks which average degree assumptions force a finite graph~$G$ to contain a given desired substructure. When this substructure is not a subgraph, but a given minor or topological minor~$H$, its presence can be forced by making the average degree of $G$ large enough in terms of~$H$ only, independently of the order of~$G$.

When $G$ is infinite, this will no longer work: since infinite trees can have arbitrarily large degrees, no minimum degree assumption can force an infinite graph to contain even a cycle. However, just as a finite tree with internal vertices of large degree has many leaves, an infinite tree of large minimum degree has many ends. The question we pursue in this paper is which notion of `degree' for ends might imply that an infinite graph whose vertices {\sl and\/} ends have large degree contains a given finite minor. (The ends of a tree should then have degree~1, so that trees are no longer counterexamples.)

Various notions of end degrees already exist, but only one of them, due to Stein\th [\the\refMayaBanff], can achieve this aim. While this was a crucial step forward, the exact notion of the `relative' end degrees she proposed still leaves room for improvement: it can be hard to verify that a given graph has large end degrees in this sense, indeed there seem to be only few sparse graphs that do.

Our aim in this paper is to develop the notion of relative end degree further, so that more graphs satisfy the premise that both their vertices and their ends have degree at least some given $k\in\N$, while keeping the notion strong enough that large vertex and end degrees can force any desired finite minor.

To indicate the notion of end degree that we have in mind, let us look at a locally finite connected graph $G$ of minimum degree~$k$ which, however, has no finite subgraph of minimum degree at least~$k$. (Such a subgraph would be good enough for us: recall that for every finite graph $H$ there exists an integer~$k$ such that every finite graph of minimum degree at least~$k$ contains $H$ as a minor.) If $G$ has a finite set $S$ of vertices such that every component $C$ of $G-S$ is such that all its vertices with a neighbour in~$S$ have degree at least~$k$ in the subgraph $G[S,C]$ of $G$ they induce together with~$S$, then $S$ and its neighbours in~$G$ induce a subgraph of minimum degree at least~$k$ in~$G$, as desired. If not, then a standard compactness argument will give us a sequence $C_0\supne C_1\supne\dots$ of `bad' components $C_n$ of $G-S_n$ for finite vertex sets~$S_n$ which identifies a unique end of~$G$. The fact that all these $C_n$ are `bad', i.e., contain a neighbour of~$S_n$ whose degree in $G[S_n,C_n]$ is less than~$k$, should then imply that this end has degree~$<k$: then the assumption that all vertices and ends of $G$ have degree at least~$k$ will imply that $G$ has a finite subgraph of minimum degree at least~$k$.

This paper is organized as follows. In Section~1 we make the above ideas precise and formulate two notions of end degree, one as indicated above, the other based on the average rather than minimum degrees of the graphs~$G[S,C]$. We shall see how the assumption of large vertex and end degrees in this sense can force finite subgraphs of given minimum or average degree. We shall also see, however, that this is rather a strong assumption: it will be nontrivial to verify that the ends of a given graph all have large degree in this sense, indeed there may not be many locally finite graphs whose vertex and end degrees are at least~$k$.%
   \COMMENT{}
   We shall therefore amend the definition of end degree so that the assumption that these degrees are large becomes weaker, thereby allowing more graphs to have large vertex and end degrees.

In Section~2 we show that for locally finite graphs with countably many ends the new notion is still strong enough to force finite subgraphs of given minimum degree. In Section~3, however, we construct an example of a graph (with continuum many ends) whose vertex and end degrees can be made arbitrarily large without forcing it to contain a finite subgraph of minimum degree~$>2$. In Section~4 we use topological methods to prove a positive result for graphs with uncountably many ends, which can also be used to strenghten our result from Section~2.

All graphs in this paper are locally finite. Notation and terminology for infinite graphs, including their ends, can be found in\th [\the\refBook].

\beginsection 1. End degrees, and statements of results

The simplest, and historically first, notion of  degree for an end is the maximum number of either disjoint or edge-disjoint rays in its graph that represent it. This notion was introduced by Halin\th [\the\refHalinEnd], who also showed that these maxima exist. Replacing in this definition the rays in an end with arcs converging to it in the Freudenthal compactification of the graph yields a more refined topological notion of end degree\th [\the\refBook,\th \the\refTopSurvey].

An assumption of large end degrees in this sense, as well as of large vertex degrees, can force some substructures in locally finite infinite graphs, such as an unspecified highly connected subgraph\th [\the\refTopSurvey,\th \the\refMayaEndDeg,\th \the\refMayaBanff,\th \the\refMayaZamora]. However, it will not yield the kind of result we are seeking here: for every integer~$k$ there are planar graphs with all vertex and end degrees at least~$k$, and hence no such degree assumption can force even a $K_5$ minor. (Indeed, if we draw the $k$-regular infinite tree in the plane and then add edges to turn each level into a cycle, we obtain a planar graph with only one end, and this end has infinite degree in all the senses discussed above.)

It was in this situation that Stein\th [\the\refMayaBanff] suggested a notion of a `relative' end degree, which would measure for a given end the ratios between the edge and the vertex boundaries of smaller and smaller regions of the graph around this end. Before we get to our main results, let us define a simplified and slightly optimized%
   \COMMENT{}
   version of this formally, and state the theorem it leads to.

Given a subgraph or vertex set $U$ in a locally finite graph~$G$, let us call
 $$V^+(U) := N(G-U)\ \ (\sub U)$$
 its {\it vertex boundary\/}. If $U$ is finite, we write
 $$d_G(U) := \sum_{v\in U} d_G(v)\> \big/\> |U|$$\noskip\medskip\noindent
 for its `average degree in~$G$'.

A {\it region\/} of $G$ is a connected induced subgraph $C$ with a finite vertex boundary. Since $G$ is locally finite, the neighbourhood $S=N(C)$ of $C$ in $G$ will then also be finite. Indeed, the regions of $G$ are precisely the subgraphs arising as components of subgraphs $G-S$ with $S\sub V(G)$ finite. We then write $G[S,C]:= G[S\cup V^+(C)]$, which is a finite graph, and call
 $$\eqalignno{\delta^+_G (C) &:= \mi_{v\in V^+(C)} d_{G[S,C]}(v)\cr
   \noalign{\noindent  the {\it minimum\/} and}
   d^+_G (C) &:= d_{G[S,C]}(V^+(C))\cr}$$%
   \COMMENT{}
 the {\it average out-degree\/} of $C$ in~$G$.

A sequence $C_0\varsupsetneq C_1\varsupsetneq\dots$ of regions is a {\it defining sequence\/} for an end~$\omega$ of~$G$ if every $C_n$ contains a ray from~$\omega$ and $\bigcap_n C_n = \es$. (This latter condition ensures that $\omega$ is the only end with a ray in every~$C_n$.)%
   \COMMENT{}
   Following the ideas of Stein\th [\the\refMayaBanff], let
 $$\eqalign{\delta^-(\omega) &:= \inff_{C_0\varsupsetneq C_1\varsupsetneq\dots}\ \lim_{n\to\infty} \ \delta^+_G(C_n)\cr
   \displaystyle d^-(\omega) &:= \inff_{C_0\varsupsetneq C_1\varsupsetneq\dots}\ \lim_{n\to\infty} \ d^+_G(C_n),\cr}$$
where both infima are taken over all defining sequences $C_0\varsupsetneq C_1\varsupsetneq\dots$ of~$\omega$ such that the limit under the infimum exists (with $\infty$ allowed).
   \COMMENT{}

Thus, $\delta^-(\omega)$ is either infinite or the least integer~$k$ such that $\omega$ has a defining sequence with all minimum out-degrees equal to~$k$. It is not hard to show for either notion that there always exists a defining sequence of regions whose out-degrees converge, possibly to~$\infty$; in the case of $\delta^-(\omega)$ the infimum is attained by a constant sequence.%
   \COMMENT{}
   Indeed, if we changed the limit in either definition to a limes inferior or a limes superior, the values of $\delta^-(\omega)$ or $d^-(\omega)$ would not change.%
   \COMMENT{}

\proclaim Theorem~\xxxNaiveThm. Let $G$ be a locally finite infinite%
   \COMMENT{}
   graph, $k\in\N$, and $q\in\Q$.
 \pitem{i} If $\delta(G)\ge k$ and $\delta^-(\omega)\ge k$ for every end $\omega$ of~$G$, then $G$ has a finite subgraph of minimum degree~$\ge k$.
 \pitem{ii} If $d_G(S)\ge q$ for every large enough%
   \footnote*{This is made precise after the statement of Theorem~\xxxCtble.}
   finite set $S\sub V(G)$, and $d^-(\omega) > q$ for every end $\omega$ of~$G$, then $G$ has a finite subgraph of average degree~$> q$.%
   \COMMENT{}
 \enditem

\proof As indicated in the Introduction. Note that the definitions of $\delta^-(\omega)$, $d^-(\omega)$ and $d_G(S)$ are chosen exactly so as to make this approach work.
   \COMMENT{}
   \endproof

%Theorem~\xxxNaiveThm\th(ii) is essentially due to Stein\th [\the\refMayaBanff]; we only fine-tuned her notion of `relative' out-degree of a region to our definition of average out-degree, to get a slightly stronger version. Theorem~\xxxNaiveThm\th (i) and its proof follow the same idea.

Natural though it may seem, Theorem~\xxxNaiveThm\ has a serious snag: it is not clear how many graphs it applies to.%
   \COMMENT{}
   Indeed, given~$k$, which locally finite graphs satisfy the premise in~(i), say, that $\delta^-(\omega)\ge k$ for every end~$\omega$? It means that {\it every\/} defining sequence $C_0\varsupsetneq C_1\varsupsetneq\dots$ for~$\omega$ is such that $\delta_G^+ (C_n)\ge k$ for every~$n$. Since it is easy to construct pathological defining sequences of ends, this is rather a strong property, and easy to foil;%
   \COMMENT{}
   see\th [\the\refMayaBanff] for a discussion.%
   \COMMENT{}

Stein\th [\the\refMayaBanff] addressed this problem in her definition of relative end degrees by restricting  the defining sequences of ends allowed in the definition of~$d^-(\omega)$. But even with these amendments it remains unclear in her paper which, if any,%
   \COMMENT{}
   graphs are such that all ends have relative degree at least some given~$k$.

It is tempting, therefore, to try a more radical cure: to change the definition of end degrees considered above by replacing their infimum with a supremum. Then an end $\omega$ will satisfy $\delta(\omega)\ge k$ as soon as there {\it exists\/} a defining sequence $C_0\varsupsetneq C_1\varsupsetneq\dots$ for it in which for infinitely many (equivalently: for all) $n$ all boundary vertices of $C_n$ have degree~$\ge k$ in~$G[S_n,C_n]$. The premise of the theorem that all these end degrees are at least~$k$ thus becomes much weaker and easier to verify, and it is easy to construct a large diversity of examples of such graphs.%
   \COMMENT{}

Thus, formally, we let
 $$\eqalign{\delta^+(\omega) &:= \su_{C_0\varsupsetneq C_1\varsupsetneq\dots}\ \lim_{n\to\infty} \ \delta^+_G(C_n)\cr
   \displaystyle d^+(\omega) &:= \su_{C_0\varsupsetneq C_1\varsupsetneq\dots}\ \lim_{n\to\infty} \ d^+_G(C_n),\cr}$$
   the suprema being taken over all defining sequences $C_0\varsupsetneq C_1\varsupsetneq\dots$ of~$\omega$ such that the limit under the infimum exists (with $\infty$ allowed).%
   \COMMENT{}
   Let us call $\delta^+(\omega)$ the {\it minimum limit degree\/} of~$\omega$, and $d^+(\omega)$ its {\it average limit degree\/}.

%Note that we no longer require formally in these definitions that the graphs $G-C_n$ be connected. This requirement would now make it harder, rather than easier, for $\delta^+(\omega)$ or $d^+(\omega)$ to be large, so theorems assuming large limit end degrees in their premise will be stronger without it. We shall prove such a theorem in Section~4. 

In its full generality, the analogue of Theorem~\xxxNaiveThm\ with this much weaker premise is now too strong to be true. To see this, let us show that every end $\omega$ of the infinite $k$-branching tree $T_k$ has a defining sequence witnessing that $\delta^+(\omega)\ge k$ (while $T_k$ clearly has no finite subgraph of minimum degree~$>1$). To obtain such a sequence, let each $C_n$ be the up-closure in~$T_k$ of a vertex~$t$ together with its lower neighbour~$t^-$. Then $V^+(C_n) = \{t^-\}$, and $S_n = N(C_n)$ consists of the upper neighbours of $t^-$ other than~$t$, and its lower neighbour. The vertex $t^-$ sends $k$ edges to $S_n$ (unless $t^-$ is the root of~$T_k$, a~case we can ignore), so $\delta^+_G (C_n) = d^+_G (C_n) = k$.%
   \COMMENT{}

We shall therefore need some restrictions in order to get positive results. One such restriction is suggested by the counterexample just discussed, which relies on the fact that the neighbourhoods of the end-defining regions used are disconnected. If we ask that these neighbourhoods be connected, or just that the graphs $G-C_n$ be connected (a~slightly weaker assumption), we get our first positive result:

\proclaim Theorem \xxxCtble. Let $G$ be a locally finite infinite%
   \COMMENT{}
   graph with at most countably many ends, and let $k\in\N$ and $q\in\Q$.
 \pitem{i} If $\delta(G)\ge k$, and if every end $\omega$ of~$G$ satisfies $\delta^+(\omega)\ge k$ witnessed by a defining sequence $C_0\varsupsetneq C_1\varsupsetneq\dots$ such that $G-C_n$ is connected for all~$n$,\penalty-200\ then $G$ has a finite subgraph of minimum degree~$\ge k$.
 \pitem{ii} If $d_G(S)\ge q$ for every large enough finite set $S\sub V(G)$, and every end $\omega$ of~$G$ satisfies $d^+(\omega) > q$ witnessed by a defining sequence $C_0\varsupsetneq C_1\varsupsetneq\dots$ such that $G-C_n$ is connected for all~$n$, then $G$ has a finite subgraph of average degree~$>q$.
 \enditem

\noindent The `large enough' in (ii) can be taken with respect to size, or to mean that only sets $S$ containing some given set $S_0$ have to satisfy $d_G(S)\ge q$. The point is that allowing singleton sets $S$ would make the condition $\foral S\: d_G(S)\ge q$ stronger than $\delta(G)\ge q$, which is not the intention.%
   \COMMENT{}
   We shall prove Theorem~\xxxCtble\ in Section~2.

In Section~4 we use topological methods to strengthen Theorem~\xxxCtble\ by dropping the requirement that the graphs $G-C_n$ be connected (Corollary~\xxxCtbleTop). The proof will assume familiarity with the more elementary proof of Theorem~\xxxCtble\ given in Section~2.

\bigskip

The assumption that $G$ should have only countably many ends is a non-trivial restriction: it is not rare that theorems that are difficult in general are much easier to prove under this assumption. So I tried for some time to prove the general version (with the additional requirement that the graphs $G-C_n$ be connected, which our earlier example has shown to be necessary)~-- but ended up finding a counterexample:\looseness=-1

\proclaim Theorem \xxxUnctble.
For every  integer $k$ there exists a locally finite graph $G$ with $\delta(G)\ge k$ all whose ends $\omega$ have minimum limit degree~$\delta^+(\omega)\ge k$, witnessed by a defining sequence $C_0\varsupsetneq C_1\varsupsetneq\dots$ such that $G-C_n$ is connected for all~$n$,%
   \COMMENT{}
   but which has no finite subgraph~$H$ with $\delta(H)>2$.

\noindent Thus, Theorem~\xxxCtble\ is best possible in this sense. The counterexample proving Theorem~\xxxUnctble\ will be described in Section~3.

\bigskip

We shall finally prove a positive result for graphs with uncountably many ends. The naive extension of Theorem~\xxxCtble\th(i) to arbitrary locally finite graphs~$G$ (which is false by Theorem~\xxxUnctble) could be rephrased, without mentioning end degrees explicitly, as saying that $G$ has a finite subgraph of minimum degree at least~$k$ as soon as $\delta(G)\ge k$ and every end has a defining sequence $C_0\varsupsetneq C_1\varsupsetneq\dots$ of regions all of minimum out-degree at least~$k$. Our positive result says that this is true if these regions can be taken from one overall set of regions of $G$ that are {\it nested\/}: such that any two of them are either disjoint or such that one contains the other. Note that, unlike in Theorem~\xxxCtble, we no longer require that the graphs $G-C_n$ be connected.

It is not uncommon for a collection of regions defining ends to be nested. For example, if $T$ is a normal spanning tree of~$G$, the up-closures in~$T$ of single vertices form a nested set of regions in which every end has a defining sequence.

In the topological terminology to be introduced in Section~4, our result takes the following form:

\proclaim Theorem~\xxxPos. Let $G$ be locally finite infinite%
   \COMMENT{}
   graph with a nested set $\C$ of regions defining a basis of~$|G|$, and let $k\in\N$ and $q\in\Q$.
 \pitem{i} If $\delta(G)\ge k$ and $\delta^+_G(C)\ge k$ for every $C\in\C$, then $G$ has a finite subgraph of minimum degree~$\ge k$.
 \pitem{ii} If $d_G(S)\ge q$ for every large enough finite set $S\sub V(G)$, and ${d^+_G(C)> q}$ for every $C\in\C$, then $G$ has a finite subgraph of average degree~$>q$.%
   \COMMENT{}
 \enditem

\noindent
Theorem~\xxxPos\ will be proved in Section~4.

%\beginsection 2. Graphs with countably many ends
\beginsection 2. Elementary positive results

To make our proof of Theorem~\xxxCtble\ easy to describe we need a few more terms.

An end of $G$ {\it lives in\/} a region $C$ if each of its rays has a tail in~$C$; then it cannot also live in another region disjoint from~$C$. Given $k\in\N$, let us call a region~$C$ {\it good\/} for the purpose of the proof of Theorem~\xxxCtble\th(i) if $G-C$ is connected and $\delta^+_G (C)\ge k$, and {\it good\/} for the purpose of the proof of Theorem~\xxxCtble\th(ii) if $G-C$ is connected and $d^+_G (C) > q$. The assumptions of $\delta^+(\omega)\ge k$ in Theorem~\xxxCtble\th(i) and of $d^+(\omega) > q$ in Theorem~\xxxCtble\th(ii) thus imply that $\omega$ has a defining sequence consisting of good regions.%
   \COMMENT{}

To prove Theorem~\xxxCtble, let us restate it more formally in its part~(ii):

\proclaim Theorem \xxxCtble. Let $G$ be a locally finite infinite%
   \COMMENT{}
   graph with at most countably many ends, and let $k\in\N$ and $q\in\Q$.
 \pitem{i} If $\delta(G)\ge k$, and if every end $\omega$ of~$G$ satisfies $\delta^+(\omega)\ge k$ witnessed by a defining sequence $C_0\varsupsetneq C_1\varsupsetneq\dots$ such that $G-C_n$ is connected for all~$n$,\penalty-200\ then $G$ has a finite subgraph of minimum degree~$\ge k$.
 \pitem{ii} If there exists a finite set $S_0\sub V(G)$ such that $d_G(S)\ge q$ for every finite set $S\sub V(G)$ containing~$S_0$, and every end $\omega$ of~$G$ satisfies $d^+(\omega) > q$ witnessed by a defining sequence $C_0\varsupsetneq C_1\varsupsetneq\dots$ such that $G-C_n$ is connected for all~$n$, then $G$ has a finite subgraph of average degree~$>q$.
   \enditem

\proof (i)
We may clearly assume that $G$ is connected. Let $\omega_1, \omega_2,\dots$ be an enumeration of the ends of~$G$. Our first aim is to select an increasing sequence $S_0\sub S_1\sub\dots$ of non-empty%
   \COMMENT{}
   finite sets of vertices satisfying the following four conditions for~$n>0$:
$$\displaylines{\hfill\hbox{\sl The end $\omega_n$ lives in a good component $C_n$ of $G-S_n$;}\hfill\llap(1)\cr
   \hfill\hbox{\sl $G[S_n]$ is connected;}\hfill\llap(2)\cr
   \hfill\hbox{\sl every component of $G-S_n$ is infinite;}\hfill\llap(3)\cr
   \hfill\hbox{\sl any good component of $G-S_{n-1}$ is also a good component of $G-S_n$.}\hfill\llap(4)}$$
As $S_0$ we take a finite set of vertices satisfying (2) and (3) for $n=0$; this can be obtained, for example, by taking a single vertex $v$ and adding all finite components of $G-v$. We now assume that we have chosen $S_0\sub\dots\sub S_n$ for some~$n$ so as to satisfy (2)--(3) if $n=0$, and (1)--(4) if $n>0$.

Let us now choose~$S_{n+1}$. If the component $C$ of $G-S_n$ in which $\omega_{n+1}$ lives is good, we let $S_{n+1} := S_n$ and $C_{n+1}:= C$. If $C$ is bad, we choose a good region $C_{n+1}\sub C$ of~$G$%
   \COMMENT{}
   in which~$\omega_{n+1}$ lives; this exists by our assumption that $\delta^+(\omega_{n+1})\ge k$.%
   \COMMENT{}
   Note that every $S_n$--$\,C_{n+1}$ path in $G$ meets~$N(C_{n+1})$ before it reaches~$C_{n+1}$.%
   \COMMENT{}
    Since $G-C_{n+1}$ is connected (as $C_{n+1}$ is good), and $G[S_n]$ is connected by~(2), we can make $N(C_{n+1})$ connected by adding finitely many finite paths from $G[S_n\cup V(C)] -C_{n+1}$. We can therefore find a finite connected subset $S'_n$ of $S_n\cup V(C)$ that contains $S_n\cup N(C_{n+1})$ but does not meet~$C_{n+1}$, and let $S_{n+1}$ be obtained from $S'_n$ by adding any finite components of $C - S'_n$. Then $C_{n+1}$ is a component of~$G-S_{n+1}$. The set $S_{n+1}$~contains~$S_n$ but is still finite, and it satisfies (1)--(4) for~$n+1$. This completes the choice of $S_0\sub S_1\sub\dots\>$.

Let us show that the above construction in fact breaks off after finitely many steps, i.e., that the sequence of sets~$S_n$ becomes stationary. If not, then $S := \bigcup_{n=0}^\infty S_n$ is infinite and spans a connected subgraph of~$G$, by~(2). As $G$ is locally finite, this subgraph contains a ray, from the end~$\omega_k$ of~$G$, say. By~(1), the end $\omega_k$ lives in the good component $C_k$ of~$S_k$, which by~(4) is still a component of $G-S$. Hence $S_k$ is a finite set of vertices separating $C_k$ from~$S$, both of which contain a ray from~$\omega_k$, a contradiction.

We have shown that our construction of sets $S_0\sub S_1\sub\dots$ breaks off after finitely many steps, with a set $S_n$ say. By~(3), every component $C$ of $G-S_n$ is infinite; let us show that it is a good region of~$G$. Since $C$ is locally finite it contains a ray. Hence some end of $G$ lives in~$C$; let $\omega_k$ be such an end with $k$ minimum. By~(1), the end~$\omega_k$ also lives in the good component $C_k$ of~$G-S_k$. If $k\le n$, then $C_k$ is still a component of~$G-S_n$, by~(4), so $C=C_k$ is good. If $n < k$, then by the maximality of~$n$ we have $S_k = S_n$. So $C_k$ is also a component of~$G-S_n$, and hence equal to~$C$.

Since every component of $G-S_n$ is good, every vertex of the finite non-empty%
   \COMMENT{}
   graph $H := G[S_n\cup N(S_n)]$ that is not in~$S_n$ has degree at least~$k$ in~$H$. But so do the vertices in~$S_n$, since they have degree~$\ge k$ in~$G$ and $H$ contains all their neighbours. Thus, having found $H$ we have completed the proof of~(i).

(ii) The proof follows the same lines as above, except that we start with $S_0$ as specified in the statement of~(ii).%
   \COMMENT{}
   Then at the end it takes a short argument to see that if every component of $G-S_n$ is good then the finite graph $H = G[S_n\cup N(S_n)]$ has average degree~$>q$. Checking this, however, is immediate from the definitions of $d_G(S_n)$ (which is at least~$q$ by assumption) and of~$d^+_G (C)$ for the components $C$ of~$G-S_n$ (which are all~$>q$), and we leave it to the reader.%
   \COMMENT{}
   \endproof

%\beginsection 3. A locally finite counterexample with continuum many ends
\beginsection 3. A counterexample

We already saw that Theorem~\xxxCtble\ does not extend to graphs with uncountably many ends if, at the same time, we drop the additional requirement on our defining sequences $C_0\supsetneq C_1\supsetneq\dots$ of ends that also the complements $G-C_n$ must be connected. However, the counterexample we saw~-- the $k$-branching tree $T_k$ with some specially chosen defining sequences for its ends~-- relied heavily on the fact that the graphs $G-C_n$ were allowed to be disconnected: this enabled us to give the $C_n$ a small vertex boundary sending many edges out: to many components of $G-C_n$, and thus not forcing high degrees outside~$C_n$. (It will become clearer below why this would otherwise likely be the case.)

We now construct a counterexample in which the graphs $G-C_n$ are connected:

\proclaim Theorem \xxxUnctble.
For every  integer $k$ there exists a locally finite graph $G$ with $\delta(G)\ge k$ all whose ends $\omega$ have minimum limit degree~$\delta^+(\omega)\ge k$, witnessed by a defining sequence $C_0\varsupsetneq C_1\varsupsetneq\dots$ such that $G-C_n$ is connected for all~$n$,%
   \COMMENT{}
   but which has no finite subgraph~$H$ with $\delta(H)>2$.

\proof
Let $k$ be given, without loss of generality $k\ge 1$. Let $T$ be the rooted tree in which every vertex $t$ has $k+1$ successors, which we think of as lying {\it above\/}~$t$. Let us call the vertices of $T$ {\it tree-vertices\/}, and its edges {\it vertical\/}. For each $t\in T$ furnish its set of $k+1$ successors with {\it horizontal\/} edges to turn it into a complete graph~$K_{k+1}$, then subdivide each of these horizontal edges once. We continue to call the subdivided edges {\it horizontal\/}, and call the new vertices the {\it subdividing vertices\/}. Call this graph~$T^+$.

Now iterate the following construction step $\omega$ times: with every subdividing vertex $s$ of the current graph identify the root of a new copy of~$T^+$, putting it above~$s$. (Keep all labels `horizontal' or `vertical' of the old graph and of the copies of $T^+$ added to it.) We will show that the resulting graph~$G$ proves the theorem.

Clearly, $\delta(G)\ge k$. To show that every end $\omega$ has a minimum limit degree as stated, let us define a set $\C$ of regions $C$ with $\delta^+_G(C) = k$ and $G-C$ connected, and such that every end has a defining sequence of regions in~$\C$. The set $\C$ will consist of one region $C_s$ for every subdividing vertex $s$ of (any of the copies of $T^+$ in)~$G$, defined as follows. Let $tt'$ be the edge of the $K_{k+1}$ of which $s$ is a subdividing vertex. Then let $C_s$ be the up-closure of $\{t,s,t'\}$ in~$G$, together with the edges $ts$ and~$st'$. This is a connected subgraph of~$G$, and $G-C_s$ is connected too. The vertex boundary of $C_s$ in $G$ consists of $t$ and~$t'$, each of which is incident with one vertical edges leaving~$C_s$ (to the common predecessor of $t$ and~$t'$ in the copy of $T^+$ containing them) and $k-1$ horizontal such edges.

Now consider an end $\omega$ of~$G$. It is represented by a ray $R$ whose vertical edges all go upwards.%
   \COMMENT{}
   For every subdividing vertex $s\in R$, the tail $sR$ of $R$ lies in~$C_s$. For every vertical edge $tt^+$ with $t$ not a subdividing vertex,%
   \COMMENT{}
   the tail $tR$ of $R$ lies in every $C_s$ whose $s$ is joined to $t$ by a horizontal edge. (There are $k\ge 1$ such vertices~$s$.) Hence every tail of $R$ lies in some~$C_s$, and these $C_s$ form a defining sequence of~$\omega$.

It remains to show that every finite subgraph $H$ of $G$ has minimum degree at most~2. Consider a highest vertex~$v$ of~$H$. If $v$ is incident with a horizontal edge $vw\in H$, then either $v$ or $w$ has degree at most~$2$ in~$H$. If not, then $v$ has degree~1 in~$H$.
   \endproof

The graph constructed in the proof of Theorem~\xxxUnctble\ can also serve as a counterexample to the naive generalization of Theorem~\xxxCtble\th(ii) to graphs with uncountably many ends: it should be easy to show that, regardless of the value of~$k$, every finite subgraph of $G$ has average degree at most~$4$, say. (This is easy to see for any complete graph in which every edge has been subdivided once, and we skip the messy details of extending the calculation to arbitrary finite subgraphs of~$G$.)%
   \COMMENT{}

%\beginsection 4. Positive results for continuum many ends
\beginsection 4. Positive results using topology

In order to prove Theorem~\xxxPos, we have to view graphs with ends in a topological setting. There is a natural topology on a locally finite graph $G$ with ends that makes it into a compact space extending the 1-complex $G$. This space~$|G|$, the {\it Freudenthal compactification\/} of~$G$, is explained in\th [\the\refBook]. All we need here is that every region $C$ of~$G$, together with the ends that live in it and the inner points of edges leaving it, forms an open subset of~$|G|$. We shall denote this open set as~$\hat C$.

Every locally finite graph $G$ has a nested set $\C$ of regions $C$ whose corresponding open sets $\hat C\sub |G|$ form a basis of~$|G|$ together with the local open stars around vertices and the open intervals of edges; for example, take the regions induced in $G$ by the up-closures of single vertices in a normal spanning tree. We then say that $\C$ {\it defines\/} this basis of~$|G|$.

\proclaim Theorem~\xxxPos. Let $G$ be locally finite infinite%
   \COMMENT{}
   graph with a nested set $\C$ of regions defining a basis of~$|G|$, and let $k\in\N$ and $q\in\Q$.
 \pitem{i} If $\delta(G)\ge k$ and $\delta^+_G(C)\ge k$ for every $C\in\C$, then $G$ has a finite subgraph of minimum degree~$\ge k$.
 \pitem{ii} If $d_G(S)\ge q$ for every large enough finite set $S\sub V(G)$, and ${d^+_G(C) > q}$ for every $C\in\C$, then $G$ has a finite subgraph of average degree~$>q$.%
   \COMMENT{}
 \enditem

\proof
We prove (i); the proof of (ii) is similar, with the same adjustments as in the proof of Theorem~\xxxCtble. 

Let $v$ be a fixed vertex of~$G$, and let $\C'$ be the subset of $\C$ consisting of the sets in $\C$ not containing~$v$. For every end~$\omega$, let $C_\omega$ be the component of $G-v$ in which $\omega$ lives. Since $\C$ defines a basis of~$|G|$ and $\hat C_\omega$ is an open subset of~$\Omega$ containing~$\omega$, there is inside~$\hat C_\omega$ a set $\hat C\owns\omega$ with $C\in\C$, and hence $C\in\C'$. So the sets~$\hat C\cap\Omega$ with $C\in\C'$ form an open cover of~$\Omega$.

Since $\Omega$ is a closed subspace of~$|G|$, and hence compact,%
   \COMMENT{}
   we can select from~$\C'$ a finite subset $\{C_1,\dots,C_n\}$ such that $\hat C_1\cup\dots\cup\hat C_n\supe \Omega$. Deleting from this set any $C_i$ contained in another~$C_j$, we may assume that $C_1,\dots,C_n$ are disjoint (using that $\C$ is nested).%
   \COMMENT{}

Since $\hat C_1,\dots,\hat C_n$ are open in~$|G|$, the subspace $X:= |G|\sm (\hat C_1\cup \dots\cup \hat C_n)$ is closed, and hence compact, but contains no end (by the choice of $C_1,\dots,C_n$). As distinct vertices of $G$ can converge only to ends in~$|G|$, this means that $X$ is a finite subgraph of~$G$. (It clearly contains every edge, including its endvertices, of which it contains an inner point.) Since~$v\in X$, it is non-empty.

As $\delta_G^+(C_i)\ge k$ for $i=1,\dots,n$, and the vertices of $X$ have at least $k$ neighbours in~$G$, the finite subgraph $H$ of $G$ spanned by $X$ and the vertex boundaries of $C_1,\dots,C_n$ has minimum degree at least~$k$. (We remark that, unlike in the proof of Theorem~\xxxCtble, $G$~may have $C_i$--$\,C_j$ edges for $i\ne j$, so the vertex boundaries of the $C_i$ need not lie in~$N(X)$. But any such edges will put their endvertices in $V^+ (C_i)$ and~$V^+ (C_j)$, so they will be edges of~$H$. The conclusion that $\delta(H)\ge k$, therefore, is still correct.)
   \endproof

Note that in the proof of Theorem~\xxxPos\ we did not use the full assumption that $\C$ defines a basis of~$|G|$, only that there exists a vertex $v\in G$ such that $\Omega$ can be covered by disjoint sets $\hat C$ not containing~$v$ such that $C$ is a good region. It would be possible, of course, to rephrase the theorem in this way.

It is also instructive to analyse our counterexample $G$ from the proof of Theorem~\xxxUnctble\ in view of Theorem~\xxxPos. By Theorem~\xxxPos, the set $\C$ of regions $C_s$ used in the example for the defining sequences of ends cannot be nested. And indeed, for every horizontal path $tst's't''$ in~$G$ the regions $C_s$ and $C_{s'}$ intersect in the entire up-closure of~$t'$. We could prevent this by taking as $C_s$ only the up-closure of $s$ itself, but then the vertex boundary of $C_s$ would only consist of~$s$, which has only two neighbours outside this set.

We finally use Theorem~\xxxPos\ to strengthen Theorem~\xxxCtble\ by dropping its connectedness requirement on the complements $G-C_n$ of regions used in defining sequences of ends:

\proclaim Corollary \xxxCtbleTop. Let $G$ be a locally finite infinite%
   \COMMENT{}
   graph with at most countably many ends, $k\in\N$, and $q\in\Q$.
 \pitem{i} If $\delta(G)\ge k$ and $\delta^+(\omega)\ge k$ for every end $\omega$ of~$G$, then $G$ has a finite subgraph of minimum degree~$\ge k$.
 \pitem{ii} If $d_G(S)\ge q$ for every large enough finite set $S\sub V(G)$, and $d^+(\omega) > q$ for every end $\omega$ of~$G$, then $G$ has a finite subgraph of average degree~$>q$.
 \enditem

\proof
Once more we prove only~(i), the proof of (ii) being similar.
For every end~$\omega$ of~$G$, pick a defining sequence $C_\omega^0\supsetneq C_\omega^1\supsetneq\dots$ witnessing that $\delta^+(\omega)\ge k$.
Let $\omega_1,\omega_2,\dots$ be a sequence in which every end of $G$ occurs infinitely often. For $n=1,2,\dots$ let $C_n$ be the first region in the sequence $C_{\omega_n}^0\supsetneq C_{\omega_n}^1\supsetneq\dots$ that has not been chosen as $C_i$ for any $i<n$ and satisfies $C_n\cap S_{n-1}=\es$ for $S_{n-1} := N(C_1)\cup\dots\cup N(C_{n-1})$. Then $C_n$ is a component of $G-S_n$. For all $i<n$ we have $S_i\sub S_n$, so $C_n$ is contained in a component of~$G-S_i$. As $C_i$, too, is a component of~$G-S_i$, either $C_n\sub C_i$ or $C_n\cap C_i = \es$.

The $C_1,C_2,\dots$, therefore, are distinct and nested, and the sequence contains infinitely many regions from each defining sequence $C_\omega^0\supsetneq C_\omega^1\supsetneq\dots\>$. So $C_1,C_2,\dots$ still defines a basis of~$|G|$, and we can apply Theorem~\xxxPos.
   \endproof

Let us conclude with an open problem which has an immediate bearing on the applicability of Theorem~\xxxPos, but which may be of interest in its own right:

\proclaim Problem. Find natural conditions for a given set of regions of $G$ defining a basis of~$|G|$ to contain nested subset that still defines a basis.

\beginsection Acknowledgement

I would like to thank Matthias Hamann for reading an early draft of this paper and suggesting a number of clarifications and corrections.

{\eightpoint
\beginsection References

\ref\refSpanningTrees
   R.\th Diestel, End spaces and spanning trees, \JCTB96 (2006), 846--854.

\font\ttt=cmtt8
\ref\refBook
     R.\th Diestel, {\it Graph theory\/}, 4th edition, Springer-Verlag
     2012. Electronic edition available at {\ttt http://diestel-graph-theory.com/}

\ref\refTopSurvey
   R.\th Diestel, Locally finite graphs with ends: a topological approach. ArXiv:0912.4213 (2012).

\ref\refMayaEndDeg
   M.\th Stein, Forcing highly connected subgraphs in locally finite graphs, \JGT54 (2007), 331--349.

\ref\refMayaBanff
   M.\th Stein, Extremal infinite graph theory, \DM311 (2011), 1472--1496.

\ref\refMayaZamora
   M.\th Stein and J.\th Zamora, The relative degree and large complete minors in infinite graphs, {\sl Electr.\ Notes in Discr.\ Math.\ \bf 37} (2011), 129--134.

\ref\refHalinEnd
   R.{\th}Halin, \"Uber die Maximalzahl frem\-der unendlicher Wege, \MN30 (1965), 63--85.

}

\bigskip\ninepoint\obeylines\parindent=0pt
Mathematisches Seminar\hfill Version 24.09.2012
Universit\"at Hamburg
Bundesstra\ss e 55
D - 20146 Hamburg
Germany

\bye